\input vanilla.sty
\pageheight {8.4in}
\pagewidth {6.9in}
\TagsOnRight
\pageno=1
\font\id=cmbx12  

\def\n{\noindent}

\def\q{\quad}
\def\m{\medpagebreak}\def\br{\bigpagebreak}

\def\I{\varphi}
\def\f{\frac}

\def\vs{\varepsilon}
\def\E{\text{\bf E}}

\def\p{\partial}

\def\br{\bigpagebreak}

\define\om{\omega}

\def\p{\partial}\define\pa{\prime}

\font\bx=cmbx10

\define\ovec{\overrightarrow}
\define\Om{\Omega}
\define\ov{\overline}
\define\CC{{\cal C}}

\def\SS{\text{\bf S}}

\define\R{\text{\bf I\!\bf R}}

\define\ff#1{\f{\p^2f}{\p{#1}^2}}
\define\\E{\text{\bf I\!\bf E}}

\def\CB{\Cal B}

\font\dubi=cmssi12
\define\BBB{{\bold B}}

\vskip20mm

\n{\id RECONSTRUCTION OF CONVEX BODIES}

\n{\id  FROM PROJECTION CURVATURE RADIUS FUNCTION}

\vskip8mm

\n{\bf R. H. Aramyan}

\font\bx=cmbx10

\vskip8mm

\n{\lineskip11.381102pt\bx In this article we pose the problem
of existence and uniqueness of convex body
for which the projection curvature radius function
coincides with given function.
We find a necessary and sufficient condition that ensures a positive
answer to both questions and
suggest an algorithm of construction of the body. Also we find
a representation of the support function of a convex body by
projection curvature radii.}

\vskip18mm

\baselineskip22.762pt

\n{\bf\S 1. INTRODUCTION}

\n Let $F(\om)$ be a function
defined on the sphere $\SS^2$. The existence and
uniqueness of convex body $\BBB\subset\R^3$
for which the mean curvature
radius at a point on $\p\BBB$ with outer normal direction $\om$
coincides with given $F(\om)$ was posed
by Christoffel (see [2],[7]). Let $R_1(\om)$ and $R_2(\om)$ be
the principal radii
of curvatures of the surface of the body at the point
with normal $\om\in\SS^2$.
Christoffel problem asked about the existence of $\BBB$
for which
$$R_1(\om)+R_2(\om)=F(\om).\tag1.1$$
The corresponding problem for
Gauss curvature $R_1(\om)R_2(\om)=F(\om)$
was posed and solved by Minkovski.
W. Blashke reduced the Christoffel problem to a
partial
differential equation of second order for the
support function (see [7]).
A. D. Aleksandrov and A. V. Pogorelov
generelized these problems, and proved the
existence and
uniqueness of convex body for which
$$G(R_1(\om),R_2(\om))=F(\om),\tag1.2$$
for a class of symmetric functions G (see [2], [9]).

\n In this paper we generalize the classic problem in a
different direction and pose a
similar problem for the projection curvature radii
of convex bodies (see [4]).
By $\CB$ we denote
the class of convex bodies $\BBB\subset\R^3$.
We need some notation.

\n $\SS^2$ -- the unit sphere in $\R^3$ (the space of
spatial directions),

\n $\SS_\om\subset\SS^2$ --
the great circle with pole at $\om\in\SS^2$,

\n $\BBB(\om)$ -- projection of $\BBB\in\CB$
onto the plane containing the origin in $\R^3$ and ortogonal to
$\om$.

\n $R(\om,\I)$ -- curvature radius of $\p\BBB(\om)$
at the point whose outer normal direction is
$\I\in\SS_\om$.

\n Let $F(\om,\I)$ be a nonnegative continuously differentable
function defined on $\{(\om,\I): \om\in\SS^2, \I\in\SS_\om\}$
(the space of "flags" see
[1]).
In this article we pose:

Problem 1.  existence and
uniqueness (up to a translation) of a convex body for which
$$R(\om,\I)=F(\om,\I)\,\,\,\,\,\,\,\,\,\,\,\,{\text and}\tag1.3$$

Problem 2. construction of that convex body.

\n It is well known (see [10]) that
a convex body $\BBB$ is determined uniquely by its
support function

\n $H(\Om)=\max\{<\Om,y>: y\in\BBB\}$
defined for $\Om\in\SS^2$, where $<\cdot,\cdot>$ denotes the standard
inner product in $\R^3$. Usually one extends
$H(\Om)$ to a function $H(x)$, $x\in\R^3$ using homogeneity:
$H(x)=|x|\,H(\Om)$, where
$\Om$ is the direction of $\ovec{Ox}$ ($O$ is the origin in $\R^3$).
Then the definition of convexity of $H(\Om)$ is written as
$$H(x+y)\le H(x)+H(y)\,\,\,\,
\,{\text for}\,{\text every}
\,\,\,\,x,y\in\R^3.$$

\n Below $\CC^k(\SS^2)$ denotes the space of $k$ times continuously
differentiable functions in $\SS^2$. A convex body $\BBB$ we call
$k$-smooth if $H(\Om)\in\CC^k(\SS^2)$.

\n Given a function $H(\Om)$ defined
for $\Om\in\SS^2$, by $H_\om(\I)$, $\I\in\SS_\om$
we denote the restriction of $H(\Om)$ to the circle
$\SS_\om$ for $\om\in\SS^2$.

\n Below we show that the
Problem 1. is equivalent
to the problem of existence of a function $H(\Om)$
defined on $\SS^2$
satisfing the
differential equation
$$H_\om(\I)+[H_\om(\I)]''_{\I\I}=F(\om,\I)\,\,\,\,
\,{\text for}\,{\text every}
\,\,\,\,\om\in\SS^2
\,\,\,\,\,{\text and }\,\,\,\,\, \I\in\SS_\om.\tag1.4$$
Note, that if restrictions of $H(\Om)$ satisfies (1.4), then
(the extention of) $H(\Om)$ is convex.

\m\n{\dubi Definition 1.1.} If for given $F(\om,\I)$
there exists $H(\Om)\in\CC^2(\SS^2)$ defined
on $\SS^2$ that satisfies (1.4), then $H(\Om)$ is called
a spherical solution of (1.4).

\n In (1.4), $H_\om(\I)$ is a flag function, so
we recall the basic concepts associated with flags
(in integral geometry the concept of a flag was first
systematically employed by R. V. Ambartzumian in [1]).

\n A flag is a pair $(\om, \I)$, where $\om\in\SS^2$
and $\I\in\SS_\om$. To each flag $(\om,\I)$ corresponds
a dual flag
$$
(\om,\I)\,\leftrightarrow\,(\om,\I)^*=(\Om,\phi),\tag1.5
$$
where $\Om\in\SS^2$ is the spatial direction same
as $\I\in\SS_\om$, while $\phi\in\SS_\Om$ is the
direction same as $\om$.
Given a flag function  $g(\om,\I)$, we denote
by $g^*$ the image of $g$ defined by
$$
g^*(\Om,\phi)=g(\om,\I),\tag1.6
$$
where $(\om,\I)^*=(\Om,\phi)$.

\m\n{\dubi Definition 1.2.} For every $\om\in\SS^2$,
(1.4) reduces to a differential equation on the circle
$\SS_\om$. Any continuous function
$G(\om,\I)$ that is a solution of (1.4)
for every $\om\in\SS^2$
we call {\dubi a flag solution}.

\m\n{\dubi Definition 1.3.} If a flag solution
$G(\om,\I)$ satisfies
$$
G^*(\Om,\phi)\,=\,G^*(\Om)\tag1.7
$$
(no dependence on the variable $\phi$), then
$G(\om,\I)$ is called {\dubi a consistent flag
solution}.

\n There is an important principle: {\dubi each
consistent flag solution $G(\om,\I)$ of (1.4)
produces a spherical solution of (1.4) via the map
$$
G(\om,\I)\to G^*(\Om,\phi)=G^*(\Om)=H(\Om),\tag1.8
$$
and vice versa:
restrictions of any spherical solution of (1.4) onto the
great circles is a consistent flag solution.}

\n Hence the problem of finding the spherical
solutions reduces to finding the
consistent flag solutions.

\n To solve the latter problem, the present paper applies
{\dubi the consistency method} first used
in [3] and [5] in an integral equations context.

\n We denote:

\n $e[\Om,\phi]$ -- the plane containing the origin of $\R^3$,
direction $\Om\in\SS^2$ and $\phi\in\SS_\Om$ ($\phi$
determine rotation of the plane around $\Om$),

\n $\BBB[\Om,\phi]$ -- projection of $\BBB\in\CB$
onto the plane $e[\Om,\phi]$,

\n $R^*(\Om,\phi)$ -- curvature radius of $\p\BBB[\Om,\phi]$
at the point whose outer normal direction is
$\Om$.

\n It is easy to see that
$$R^*(\Om,\phi)\,=\,R(\om,\I),$$
where $(\Om,\phi)$ is the flag dual to $(\om,\I)$.

\n Note, that in the Problem 1. uniquness (up to a translation) follows from
the classical uniqueness result on
Christoffel problem, since
$$R_1(\Om)+R_2(\Om)=\f1\pi\int_0^{2\pi}R^*(\Om,\phi)\,d\phi.\tag1.9$$

\n In case $F(\om,\I)\ge0$ is nonnegative, the
equation (1.4) has the following geometrical
interpretation.

\n It follows from [4] that homogeneus fonction $H(x)=|x|H(\Om)$,
where $H(\Om)\in \CC^2(\SS^2)$, is
convex if and only if
$$H_\om(\I)+[H_\om(\I)]''_{\I\I}\ge0\,\,\,\,\,{\text for}\,{\text  every}
\,\,\,\,\,\om\in\SS^2
\,\,\,\,\,{\text and }\,\,\,\,\, \I\in\SS_\om,\tag1.10$$
where $H_\om(\I)$ is  the restriction of $H(\Om)$ onto $\SS_\om$.

\n So in case $F(\om,\I)\ge0$, it follows from (1.10), that
if $H(\Om)$ is a
spherical solution of (1.4) then
its homogeneus extention $H(x)=|x|\,H(\Om)$ is convex.

\n It is well known from convexity theory that if a function $H(x)$
is convex then there is a unique convex body
$\BBB\subset\R^3$
with support function $H(x)$ and $F(\om,\I)$ is the projection
curvature radius function of $\BBB$ (see [8]).

\n The support function
of each parallel shifts (translation)
of that body $\BBB$ will again be a spherical solution of (1.4).
By uniqueness, every two spherical solutions of (1.4) differ
by a summand $<a,\Om>$, where $a\in\R^3$.
Thus we proved the following theorem.

\m\n{\dubi Theorem 1.1. Let $F(\om,\I)\ge0$ be a nonnegative function
defined on $\{(\om,\I): \om\in\SS^2, \I\in\SS_\om\}$.
If the equation (1.4)
has a spherical solution $H(\Om)$ then there exists a convex
body $\BBB$ with projection curvature radius function $F(\om,\I)$,
whose support function is $H(\Om)$.
Every spherical
solutions of (1.4) has the form
$H(\Om)+<a,\Om>$, where $a\in\R^3$,
each being the support function of a translation
of the convex body $\BBB$ by $\ovec{Oa}$.
}

\n The converse statement is also true.
It follows from the theory for $2$-dimension (see [8]), that the
support function $H(\Om)$ of a 2-smooth convex body $\BBB$ satisfies (1.4)
for $F(\om,\I)=R(\om,\I)$, where $R(\om,\I)$ is the projection
curvature radius function of $\BBB$.

\n Before going to the main result, we make some
remarks. The purpose of the present paper is to find
a necessary and sufficient condition that ensures a positive
answer to both Problems 1,2 and
suggest an algorithm of construction of the body $\BBB$
by finding a representation of the support function
in terms of projection curvature radius function.
This happens to be a spherical solution of the equation (1.4).
In this paper the support function of a convex body $\BBB$
is considered with respect to a special choice of the origin
$O^*$.
It turns out
that each $1$-smooth convex body $\BBB$ has the special point $O^*$
we will call the centroid of $\BBB$ (see Theorem 6.1).
The centroid coincides with the centre of
symmetry for centrally symmetrical
convex bodies. 

\n For convex bodies $\BBB$ with positive Gaussian curvature one can
define the centroid as follows: within $\BBB$
there exists a unique point $O^*$ such that (see Lemma 6.1)
$$\int_{\SS_\Om}<\ovec{O^* P_\Om(\tau)},\Om>\,d\tau=0 \,\,\,\,
{\text for }\,{\text every}\,\,\,\,\,\Om\in\SS^2,$$
where 
$P_\Om(\tau)$ is the point on
$\p\BBB$ whose outer normal has the direction $\tau\in\SS_\Om$,
$d\tau$ is the usual angular measure on $\SS_\Om)$. The set of points
$\{P_\Om(\tau), \tau\in\SS_\Om\}$ we will call the belt of $\BBB$
with normal $\Om$.

\n Throughout the paper (in particular, in Theorem
1.2 that follows) we use usual spherical coordinates
$\nu,\tau$ for points $\om\in\SS^2$ based on a
choice of a North Pole $\Cal N\in\SS^2$ and
a reference point $\tau=0$ on the
equator $\SS_{\Cal N}$. 
We put $\nu=\f{\pi}2-
(\widehat{\om,\Cal N})$ so that the points
$(0,\tau)$ lie on the equator $\SS_{\Cal N}$.
The point with coordinates $\nu,\tau$
we will denote by $(\nu,\tau)_{\Cal N}$.
On each $\SS_\om$ we choose $E$=the
direction East for the reference point
and the anticlockwise direction as positive.

\n Now we describe the main result.

\m\n{\dubi Theorem 1.2. The support function of any
$3$-smooth convex body $\BBB$ with respect to the centroid
$O^*$ has the representation
$$H(\Om)=\f1{4\pi}\int_0^{2\pi}\left[\int_0^{\f\pi2}R((0,\tau)_\Om,\I)\,
\cos\I\,d\I\right]\,d\tau+\f1{8\pi^2}\int_0^{2\pi}
\left[\int_{-\f\pi2}^{\f\pi2}R((0,\tau)_\Om,\I)\,
((\pi+2\I)\cos\I-2\sin^3\I)\,d\I\right]\,d\tau-$$
$$-\f1{2\pi^2}\int_0^{\f\pi2}\f{sin \nu}{\cos^2\nu}d\nu\int_0^{2\pi}d\tau
\int_0^{2\pi}R((\nu,\tau)_\Om,\I)\sin^3{\I}\,d\I\tag1.11$$
where $R(\om,\I)$ is the projection curvature radius function
of $\BBB$, on $\SS_\om$ we measure $\I$ from the East
direction with respect to $\Om$.
(1.11) is a spherical solution
of the equation (1.4) for $F(\om,\I)=R(\om,\I)$.}

\n Remark, that the order of integration in the last integral of (1.11)
is important.

\n Obviously Theorem 1.2 suggests a practical algorithm
of reconstruction of convex bodies from
projection curvature radius function $R(\om,\I)$
by calculation of the support function $H(\Om)$.

\n We turn to Problem 1.
Let $R(\om,\I)$ be the projection curvature radius function
of a convex body
$\BBB$. Then $F(\om,\I)\equiv R(\om,\I)$ necessarily satisfies the following
conditions:

\n 1. 
$$\int_0^{2\pi}F(\om,\I)\,\sin\I\,d\I=
\int_0^{2\pi}F(\om,\I)\,\cos\I\,d\I=0,\tag1.12$$
for every $\om\in\SS^2$ and any reference point on $\SS_\om$
(follows from equation (1.4), see also in [8]).

\n 2. For every direction $\Om\in\SS^2$
$$\int_0^{2\pi}[F^*((\nu,\tau)_\Om,N)]'_{\nu=0}
\,d\tau=0,\tag1.13$$
where $F^*(\Om,\phi)=F(\om,\I)$ (see (1.6)) and
$N=E+\f\pi2$ is the North direction at the
point $(\nu,\tau)_\Om$ with respect $\Om$ (Theorem 5.1).

\n Let $F(\om,\I)$ be a nonnegative continuously differentable
function defined on $\{(\om,\I): \om\in\SS^2, \I\in\SS_\om\}$.

\n Using (1.11), we construct a function
$\ov F(\Om)$ defined on $\SS^2$:
$$\ov F(\Om)=\f1{4\pi}\int_0^{2\pi}
\left[\int_0^{\f\pi2}F((0,\tau)_\Om,\I)\,
\cos\I\,d\I\right]\,d\tau+\f1{8\pi^2}\int_0^{2\pi}
\left[\int_{-\f\pi2}^{\f\pi2}F((0,\tau)_\Om,\I)\,
((\pi+2\I)\cos\I-2\sin^3\I)\,d\I\right]\,d\tau-$$
$$-\f1{2\pi^2}\int_0^{\f\pi2}
\f{sin \nu}{\cos^2\nu}d\nu\int_0^{2\pi}d\tau
\int_0^{2\pi}F((\nu,\tau)_\Om,\I)
\sin^3{\I}\,d\I\tag1.14$$
Note that the last integral converges if the condition (1.13) is satisfied
(see (5.7) and (5.8)).

\proclaim{Theorem 1.3} A nonnegative
continuously differentable
function $F(\om,\I)$ defined on $\{(\om,\I): \om\in\SS^2, \I\in\SS_\om\}$
represents the
projection curvature radius function of some convex body
if and only if $F(\om,\I)$
satisfies the conditions (1.12), (1.13)
and
$$\ov F_\om(\I)+[\ov F_\om(\I)]''_{\I\I}=F(\om,\I)\,\,\,\,
\,{\text for}\,{\text every}
\,\,\,\,\om\in\SS^2
\,\,\,\,\,{\text and }\,\,\,\,\, \I\in\SS_\om.\tag1.15$$
where $\ov F_\om(\I)$ is the restriction of $\ov F(\Om)$ (given by (1.14))
onto $\SS_\om$.
\endproclaim

\n Note that, in [6] the same problem for centrally symmetrical convex
bodies was posed and a necessary and sufficient condition
ensuring a positive answer found.

\m\n{\bf \S 2. GENERAL FLAG SOLUTION OF (1.4)}

\n We fix $\om\in\SS^2$ and a pole $\Cal N\in\SS^2$ and try
to solve (1.4) as a
differential equation of second order on the circle $\SS_\om$.

\n We start with two results from [8].

1. For any smooth convex domain
$D$ in the plane
$$h(\I)=\int_0^\I R(\psi)\,\sin(\I-\psi)\,d\psi,\tag2.1$$
where $h(\I)$ is the support function of $D$ with respect
to a point $s\in\p D$. In (2.1) we measure $\I$
from the normal direction
at $s$, $R(\psi)$ is the curvature radius of $\p D$ at the point
with normal $\psi$.

2. (2.1) is a solution of the following differential equation
$$R(\I)=h(\I)+h''(\I).\tag2.2$$

\n One can easy verify that (also it follows from (2.2) and (2.1))
$$G(\om,\I)=\int_0^\I F(\om,\psi)\,\sin(\I-\psi)\,d\psi,\tag2.3$$
is a flag solution of the equation (1.4).

\m\n{\dubi Theorem 2.1. Every flag solution of (1.4)
has the form
$$
g(\om,\I)=\int_0^\I F(\om,\psi)\,\sin(\I-\psi)\,d\psi
+C(\om)\cos{\I}+ S(\om)\sin{\I}\tag2.4
$$
whera $C_n$ and $S_n$ are some real coefficients.}

\m\n{\dubi Proof:} Every continuous flag solution
of (1.4) is a sum of $G(\om,\I)+g_0(\om,\I)$, where
$g_0(\om,\I)$ is a flag solution of the corresponding
homogeneous equation
$$
H_\om(\I)+[H_\om(\I)]''_{\I\I}=0\,\,\,\,\,{\text for}\,{\text every}
\,\,\,\om\in\SS^2
\,\,\,\,\,{\text and }\, \I\in\SS_\om.\tag2.5$$
We look for the general flag solution of (2.5) as a
Fourier series
$$
g_0(\om,\I)=\sum_{n=0,1,2,...}
[C_n(\om)\cos{n\I}+S_n(\om)\sin{n\I}].\tag2.6
$$
After substitution of (2.6) into (2.5) we obtain
that $g_0(\om,\I)$ satisfy (2.5) if and only if
it has the form
$$g_0(\om,\I)=C_1(\om)\cos{\I}+ S_1(\om)\sin{\I}.$$
Theorem 2.1 is proved.

\m\n{\bf \S 3. THE CONSISTENCY CONDITION}

\n Now we consider $C=C(\om)$ and $S=S(\om)$
in (2.4) as functions of $\om=(\nu,\tau)$ and try
to find $C(\om)$ and $S(\om)$ from the condition that
$g(\om,\I)$ satisfies (1.7). We write $g(\om,\I)$ in
dual coordinates i.e. $g(\om,\I)=g^*(\Om,\phi)$
and require that $g^*(\Om,\phi)$
should not depend on $\phi$ for every $\Om\in\SS^2$,
i.e. for every $\Om\in\SS^2$
$$
(g^*(\Om,\phi))'_\phi=\left(G(\om,\I)+C(\om)\cos{\I}+
S(\om)\sin{\I}\right)'_\phi=0,\tag3.1
$$
where $G(\om,\I)$ was defined in (2.3).

\n Here and below $(\cdot)'_\phi$ denotes the
derivative corresponding to right screw rotation
around $\Om$.

\n Termwise differentiation with use of expressions
(see [5])
$$
\tau'_\phi=\f{\sin\I}{\cos\nu},\q
\I'_\phi=-\tan\nu\,\sin\I,\q
\nu'_\phi=-\cos\I,\tag3.2
$$
after a natural grouping of the summands
in (3.1), yields the Fourier series of $-(G(\om,\I))^\pa_\phi$
(a detailed derivation is contained
in [5] and [3]). By uniqueness of the Fourier coefficients

$$\cases
(C(\om))'_\nu+\dsize\f{(S(\om))'_\tau}
{\cos\nu}+\tan\nu\,C(\om)=
\f1\pi\int_0^{2\pi}A(\om,\I)\,\cos{2\I}\,d\I\\
(C(\om))'_\nu-\dsize\f{(S(\om))'_\tau}{\cos\nu}-\tan\nu\,
C(\om)=\f1{2\pi}\int_0^{2\pi}A(\om,\I)\,d\I \\
(S(\om))'_\nu-\dsize\f{(C(\om))'_\tau}
{\cos\nu}+\tan\nu\,S(\om)=
\f1\pi\int_0^{2\pi}A(\om,\I)\,\sin{2\I}\,d\I,
\endcases\tag3.3
$$
where
$$
A(\om,\I)=
\int_0^\I[F(\om,\psi)'_\phi\,\sin{(\I-\psi)}+
F(\om,\psi)\,\cos{(\I-\psi)}\,\I'_\phi]\,d\psi.\tag3.4$$

\m\n{\bf \S 4. AVERAGING}

\n Let $H$ be a spherical solution
of (1.4), i.e. restriction of $H$ onto the great
circles is a consistent flag solution of (1.4).
By Theorem 1.1 there exists a convex
body $\BBB$ with projection curvature
radius function $R(\om,\I)=F(\om,\I)$,
whose support function is $H(\Om)$.

\n To calculate $H(\Om)$ we take $\Om\in\SS^2$ for
the pole $\Om=\Cal N$. Returning to the formula
(2.4) for every $\om=(0,\tau)_\Om\in\SS_\Om$ we have
$$
H(\Om)=\int_0^{\f\pi2} R(\om,\psi)\,\sin(\f\pi2-\psi)\,d\psi
+S(\om),\tag4.1
$$
We integrate both sides of (4.1) with respect to
uniform angular measure $d\tau$ over $[0,2\pi)$ to get
$$
2\pi H(\Om)=\int_0^{2\pi}\int_0^{\f\pi2}
R((0,\tau)_\Om,\psi)\,\cos\psi\,d\psi\,d\tau
+\int_0^{2\pi}S((0,\tau)_\Om)\,d\tau.\tag4.2
$$
Now the problem is to calculate
$$
\int_0^{2\pi}S((0,\tau)_\Om)\,d\tau=\ov{S}(0).
\tag4.3
$$
We are going to integrate both sides of (3.3) and
(3.4) with respect to $d\tau$ over $[0,2\pi)$.

\n For $\om=(\nu,\tau)_\Om$, where $\nu\in[0,
\f\pi2)$ and $\tau\in(0,2\pi)$ (see (3.5))
we denote
$$
\ov S(\nu)=
\int_0^{2\pi}S((\nu,\tau)_\Om)\,d\tau,\tag4.4
$$
$$
A(\nu)=\f1{\pi}\int_0^{2\pi}\,d\tau\int_0^{2\pi}\left[
\int_0^\I[R(\om,\psi)'_\phi\,\sin{(\I-\psi)}+
R(\om,\psi)\,\cos{(\I-\psi)}\,\I'_\phi]\,
d\psi\right]\sin{2\I}\,d\I.\tag4.5$$

\n Integrating both sides of (3.3) and (3.4) and
taking into account that
$$
\int_0^{2\pi}(C(\nu,\tau)_\Om)'_\tau\,d\tau=0
$$
for $\nu\in[0,\f\pi2)$ we get
$$\ov S'(\nu)+\tan\nu\,\ov S(\nu)=A(\nu).\tag4.6$$
Thus we have differential equation (4.6) for
unknown coefficient $\ov S(\nu)$.

\m\n{\bf \S 5. BOUNDARY CONDITION FOR DIFFERENTIAL EQUATION (4.6)}

\n We have to find $\ov S(0)$ given by (4.3). It follows from (4.6) that
$$\left(\f{\ov S(\nu)}{\cos\nu}\right)'=\f{A(\nu)}{\cos\nu}.\tag5.1$$
Integrating both sides of (5.1) with respect to $d\nu$ over $[0,\f\pi2)$
we obtain
$$\ov S(0)=\left.\f{\ov S(\nu)}{\cos\nu}\right|_{\f\pi2}-\int_0^{\f\pi2}
\f{A(\nu)}{\cos\nu}\,d\nu.\tag5.2$$
Now, we are going to calculate
$\left.\f{\ov S(\nu)}{\cos\nu}\right|_{\f\pi2}$.

\n It follows from (2.4) that
$$\ov S(\nu)=\f1\pi\int_0^{2\pi}\int_0^{2\pi}\left[
H_\om(\I)-\int_0^{\I} R(\om,\psi)\,
\sin(\I-\psi)\,d\psi\right]\sin\I\,d\I\,d\tau=$$
$$=\f1\pi\int_0^{2\pi}\int_0^{2\pi}
H_\om(\I)\sin\I\,d\I\,d\tau-
\f1{2\pi}\int_0^{2\pi}\int_0^{2\pi}
R(\om,\psi)\left((2\pi-\psi)\cos\psi+\sin\psi
\right)\,d\psi\,d\tau.\tag5.3$$
Let $\I\in\SS_\om$ be the direction that
corresponds to $\I\in [0,2\pi)$, for $\om=(\nu,\tau)_\Om$.
As a point of
$\SS^2$, let $\I$ have spherical coordinates $u,t$
with respect $\Om$. By the sinus theorem of spherical geometry

$$\cos\nu\,\sin\I=\sin u.\tag5.4$$
From (5.4) we get
$$(u)'_{\nu=\f\pi2}=-\sin\I.\tag5.5$$
Using (5.5), for a fix $\tau$ we write
a Taylor expession
at a neighbourhood of the point $\nu=\f\pi2$:
$$H_{(\nu,\tau)_\Om}(\I)=H((0,\I+\tau)_\Om)+
H'_\nu((0,\I+\tau)_\Om)\,\sin\I\,(\f\pi2-\nu)+o(\f\pi2-\nu).\tag5.6$$
Similarly, for $\psi\in[0,2\pi)$ we get
$$R((\nu,\tau)_\Om,\psi)=R((\f\pi2,\tau)_\Om,\psi+\tau)+
R'_\nu((\f\pi2,\tau)_\Om,\psi+\tau)
\,\sin\psi\,(\f\pi2-\nu)+o(\f\pi2-\nu).\tag5.7$$
Substituting (5.6) and (5.7) into (5.3) and taking into account
the easy equalities
$$\int_0^{2\pi}\int_0^{2\pi}
H((0,\I+\tau)_\Om)\sin\I\,d\I\,d\tau=0$$
and
$$\int_0^{2\pi}\int_0^{2\pi}
R((\f\pi2,\tau)_\Om,\psi+\tau)
\left((2\pi-\psi)\cos\psi+\sin\psi
\right)\,d\psi\,d\tau=0\tag5.8$$
we obtain
$$\lim_{\nu\to\f\pi2}\f{\ov S(\nu)}{\cos\nu}=
\f1\pi\int_0^{2\pi}\int_0^{2\pi}
H'_\nu((0,\I+\tau)_\Om)
\,\sin^2\I\,d\I\,d\tau-$$
$$-\f1{2\pi}\int_0^{2\pi}\int_0^{2\pi}
R'_\nu((\f\pi2,\tau)_\Om,\psi+\tau)
\,\sin\psi\,
\left((2\pi-\psi)\cos\psi+\sin\psi
\right)\,d\psi\,d\tau=$$
$$=\int_0^{2\pi}
H'_\nu((0,\tau)_\Om)
\,d\tau-
\f34\int_0^{2\pi}
[R^*((\nu,\tau)_\Om,N)]'_{\nu=0}
d\tau.\tag5.9$$

\m\n{\dubi Theorem 5.1. For every $3$-smooth convex body $\BBB$
and any direction $\Om\in\SS^2$, we have
$$\int_0^{2\pi}[R^*((\nu,\tau)_\Om,N)]'_{\nu=0}
\,d\tau=0,\tag5.10$$
where ${\nu,\tau}$ is the spherical coordinates with respect to $\Om$,
where $N=E+\f\pi2$ is the North direction at the
point $(\nu,\tau)_\Om$ with respect $\Om$.}

\m\n{\dubi Proof.} 
Using spherical geometry, one
can prove that (see also (1.4))
$$[R^*((\nu,\tau)_\Om,N)]'_{\nu=0}=
\left[H((\nu,\tau)_\Om)+
H''_{\I\I}((\nu,\tau)_\Om)\right]'_{\nu=0}=$$
$$\left[H((\nu,\tau)_\Om)+H''_{\tau\tau}\f1{\cos^2\nu}-
H'_\nu\tan\nu\right]'_{\nu=0}=
[H''_{\tau\tau}]'_{\nu=0},\tag5.11$$
where $H(\Om)$ is the supporet function of $\BBB$.
After integration (5.11) we get
$$\int_0^{2\pi}
[R^*((\nu,\tau)_\Om,N)]'_{\nu=0}\,d\tau=
\int_0^{2\pi}[H''_{\tau\tau}]'_{\nu=0}\,d\tau=0.$$

\m\n{\bf \S 6. CENTROID OF A CONVEX BODY}

\n Let $\BBB$ be a convex body in $\R^3$ and
$Q\in\R^3$ be a point. By $H_Q(\Om)$ we denote the support function
of $\BBB$  with respect to $Q$.

\m\n{\dubi Theorem 6.1. For a given $1$-smooth convex body $\BBB$
there is a point $O^*\in\R^3$ such that
$$\int_0^{2\pi}\left[H_{O^*}((\nu,\tau)_\Om)\right]'_{\nu=0}\,
d\tau=0 \,\,\,\,
{\text for }\,{\text every}\,\,\,\,\,\Om\in\SS^2,\tag6.1$$
where ${\nu,\tau}$ are the spherical coordinates with respect $\Om$.
}

\m\n{\dubi Proof.} For a given $\BBB$ and a point $Q\in\R^3$
by $K_Q(\Om)$ we denote the following function defined on $\SS^2$
$$K_Q(\Om)=\int_0^{2\pi}\left[H_{Q}((\nu,\tau)_\Om)\right]'_{\nu=0}\,
d\tau.$$
$K_Q(\Om)$ is a continuous odd function with maximum $\ov K(Q)$
$$\ov K(Q)=\max_{\Om\in\SS^2}K_Q(\Om).$$
It is easy to see that $\ov K(Q)\to\infty$ for
$|Q|\to\infty$. Since $\ov K(Q)$ is a continuous
so there is a point $O^*$ for which
$$\ov K(O^*)=\min\ov K(Q).$$
Let $\Om^*$ be a (say unique) direction of maximum i.e.
$$\ov K(O^*)=\max_{\Om\in\SS^2}K_{O^*}(\Om)=K_{O^*}(\Om^*).$$
If $\ov K(O^*)=0$ the theorem is proved. For the case
$\ov K(O^*)=a>0$ let $O^{**}$ be the point for which
$\ovec{O^*O^{**}}=\vs\,\Om^*$. It is easy to understand that
$H_{O^{**}}(\Om)=H_{O^*}(\Om)-\vs(\Om,\Om^*)$, hence
for a small $\vs>0$ we find
that $\ov K(O^{**})=a-2\pi\vs$ which is contrary to definition of $O^*$.
So $\ov K(O^*)=0$. For the case where there are two or
more directions of maximum one can
apply a similar argument. The theorem is proved.

\n The point $O^*$ we will call {\dubi the centroid} of the convex body
$\BBB$. Theorem 6.2 below gives a clearer geometrical interpretation
to that concept.

\n Let
$P_\Om(\tau)$ be the point on
$\p\BBB$ whose outer normal has the direction $\tau\in\SS_\Om$.

\m\n{\dubi Lemma 6.1. For every $2$-smooth convex body $\BBB$
with positive Gaussian curvature
and any direction $\Om\in\SS^2$, we have
$$\int_0^{2\pi}\left[H_{Q}((\nu,\tau)_\Om)\right]'_{\nu=0}\,d\tau=
\int_{\SS_\Om}<\ovec{Q P_\Om(\tau)},\Om>\,d\tau,\tag6.2$$
where $Q$ is a point of $\R^3$ and
$d\tau$ is the usual angular measure on $\SS_\Om)$.}

\m\n{\dubi Proof.} Let $B[\Om,\tau]$ be the projection of $\BBB$
onto the plane $e[\Om,\tau]$ (containing $Q$ and the
directions $\Om\in\SS^2$ and $\tau\in\SS_\Om$) and $P^*(\tau)$ be
the point on $\p B[\Om,\tau]$
with outer normal $(0,\tau)_\Om$. For the support function
of $B[\Om,\tau]$ (equivalently for the restriction of $H_Q(\Om)$
onto $e[\Om,\tau]$) we have
$$[H_Q((\nu,\tau)_\Om)]'_{\nu=0}=[|\ovec{Q P^*(\tau)}|\,
\cos{(\nu-\nu_o)}+H_{P^*}(\nu)]'_{\nu=0}=|\ovec{Q P^*(\tau)}|\,
\sin{\nu_o}=<\ovec{Q P^*(\tau)},\Om>,\tag6.3$$
where $H_{P^*}(\nu)$ is the support function of $B[\Om,\tau]$
with respect to the point $P^*(\tau)\in\p B[\Om,\tau]$
and $(\nu_0,\tau)_\Om$ is the direction of
$\ovec{Q P^*(\tau)}$.
The statment
$[H_{P^*}(\nu)]'_{\nu=0}=0$ was proved in [8].
Integrating (6.3) and taking into account that
$<\ovec{Q P^*(\tau)},\Om>=<\ovec{Q P_\Om(\tau)},\Om>$ we get (6.2).

\n Theorem 6.1 and Lemma 6.1 imply the following Theorem.

\m\n{\dubi Theorem 6.2. For
a 2-smooth convex body $\BBB$ with positive Gaussian curvature
we have
$$\int_0^{2\pi}<\ovec{O^* P_\Om(\tau)},\Om>\,d\tau
=0 \,\,\,\,
{\text for }\,{\text every}\,\,\,\,\,\Om\in\SS^2,\tag6.4$$
where $O^*$ is the centroid of $\BBB$.}

\n One can consider the last statement as a definition of
the centroid of $\BBB$.

\m\n{\bf \S 7. A REPRESENTATION FOR SUPPORT FUNCTION OF CONVEX BODIES}

\n Let $O^*$ be the centroid of the convex body
$\BBB$ (see \S 6). Now we take
$O^*$ for the origin of $\R^3$.
Below $H_{O^*}(\Om)$ we will simply denote by $H(\Om)$.

\n By Theorem 6.1, Theorem 5.1 and Lemma 6.1 we have the boundary
condition (see (5.9))
$$\left.\f{\ov S(\nu)}{\cos\nu}\right|_{\f\pi2}=0.\tag7.1$$
Substituting (5.2) into (4.2) we get
$$
2\pi H(\Om)=\int_0^{2\pi}\int_0^{\f\pi2}
R((0,\tau)_\Om,\psi)\,\cos\psi\,d\psi\,d\tau
-\int_0^{\f\pi2}
\f{A(\nu)}{\cos\nu}\,d\nu=\int_0^{2\pi}\int_0^{\f\pi2}
R((0,\tau)_\Om,\psi)\,\cos\psi\,d\psi\,d\tau-$$
$$-\f1\pi\int_0^{\f\pi2}
\f{d\nu}{\cos\nu}
\int_0^{2\pi}\,d\tau\int_0^{2\pi}\left[
\int_0^\I[R(\om,\psi)'_\phi\,\sin{(\I-\psi)}+
R(\om,\psi)\,\cos{(\I-\psi)}\,\I'_\phi]\,
d\psi\right]\sin{2\I}\,d\I.\tag7.2$$
Using expressions (3.2) and integrating by $d\I$ yields
$$2\pi\,H(\Om)=\int_0^{2\pi}\int_0^{\f\pi2}
R((0,\tau)_\Om,\psi)\,\cos\psi\,d\psi\,d\tau
+\f1\pi\int_0^{\f\pi2}
\f{d\nu}{\cos\nu}
\int_0^{2\pi}\,d\tau\int_0^{2\pi}\left[
R(\om,\psi)'_\nu\,I+
R(\om,\psi)\,\tan\nu\,II\right]\,d\psi,\tag7.2$$
where
$$II=\int_\psi^{2\pi}\sin{2\I}\,\cos(\I-\psi)\,\sin\I\,d\I=
\left[\f{(2\pi-\psi)\cos\psi}4+
\f{\sin\psi(1+\sin^2\psi)}4-\sin^3\psi\right],$$
and
$$I=\int_\psi^{2\pi}\sin{2\I}\,\sin(\I-\psi)\,\cos\I\,d\I=
\left[\f{(2\pi-\psi)\cos\psi}4+
\f{\sin\psi(1+\sin^2\psi)}4\right].\tag7.3$$
Integrating by parts (7.2) we get
$$2\pi\,H(\Om)=\int_0^{2\pi}\int_0^{\f\pi2}
R((0,\tau)_\Om,\psi)\,\cos\psi\,d\psi\,d\tau
-\f1\pi\int_0^{\f\pi2}
{d\nu}
\int_0^{2\pi}\,d\tau\int_0^{2\pi}R(\om,\psi)
\f{\sin\nu\,\sin^3\psi}{\cos^2\nu}\,d\psi-$$
$$-\f1\pi\int_0^{\f\pi2}
\,d\tau\int_0^{2\pi}R((0,\tau)_\Om,\psi)
I\,d\psi+
\lim_{a\to\f\pi2}\f1{\pi\cos a}\int_0^{\f\pi2}
\,d\tau\int_0^{2\pi}R((a,\tau)_\Om,\psi)\,I
\,d\psi.\tag7.4$$
Using (5.7), Theorem 5.1 and taking into account that
$$\int_0^{2\pi}I\,d\psi=0$$
we get
$$2\pi\,H(\Om)=\int_0^{2\pi}\int_0^{\f\pi2}
R((0,\tau)_\Om,\psi)\,\cos\psi\,d\psi\,d\tau-$$
$$-\f1\pi\int_0^{\f\pi2}
{d\nu}
\int_0^{2\pi}\,d\tau\int_0^{2\pi}R(\om,\psi)
\f{\sin\nu\,\sin^3\psi}{\cos^2\nu}\,d\psi-
\f1\pi\int_0^{\f\pi2}
\,d\tau\int_0^{2\pi}R((0,\tau),\psi)
I\,d\psi.\tag7.5$$
From (7.5), by (1.12) we obtain (1.11). Theorem 1.2 is proved.

\m\n{\bf \S 8. PROOF OF THEOREM 1.3}

\m\n{\dubi Proof.} Necessity:
let $F(\om,\I)$ be the projection curvature radius function
of a convex body
$\BBB$, then it satisfies (1.12) (see [8]),
the condition (1.13) (Theorem 5.1)
and the condition (1.14) (Theorem 1.2).

\n Sufficiently: let $F(\om,\psi)$ be a
nonnegative
continuous differentable function
satisfing the conditions (1.12), (1.13), (1.15).
By means of (1.14) we construct the function $\ov F(\Om)$
defined on $\SS^2$ as in (1.14).
According to (1.15), $\ov F(\Om)$ is a convex function
hence there exists a convex body $\BBB$
with support function $\ov F(\Om)$.
The same (1.15) implies that $F(\om,\I)$ is the projection
curvature radius of $\BBB$.

\n I would like to express my gratitude to Professor R. V. Ambartzumian
for helpful remarks.

\baselineskip11.381102pt

\br\n{\bf R E F E R E N C E S}

\br\item{1. }R. V. Ambartzumian, ``Factorization
Calculus and Geometrical Probability", Cambridge
Univ. Press, Cambridge, 1990.
          
\item{2. }A. D. Alexandrov,
``Uniqueness theorems for surfaces in the large"[in Russian],
Vesti, LGU, no. 19,  1956.

\item{3. }R. H. Aramyan, ``An approach to generalized
Funk equations, I"
[in Russian],
Izv. Akad. Nauk Armenii. Matematika, [English
translation: Journal of Contemporary Math. Anal.
(Armenian Academy of Sciences)], vol. 36, no. 1,
pp. 47 -- 58, 2001.

\item{4. }R. H. Aramyan, ``Curvature radii of
planar projections of convex bodies in $R^n$" [in
Russian], Izv. Akad. Nauk Armenii. Matematika,
[English translation: Journal of Contemporary
Math. Anal. (Armenian Academy of Sciences)],
vol. 37, no. 1, pp. 2 -- 14, 2002.

\item{5. }R. H. Aramyan, ``Generalized Radon
transform with an application in convexity theory"
[in Russian],
Izv. Akad. Nauk Armenii. Matematika, [English
translation: Journal of Contemporary Math. Anal.
(Armenian Academy of Sciences)], vol. 38, no. 3,
2003.

\item{6. }R. H. Aramyan, ``Reconstruction of centrally
symmetrical convex bodies by projection curvature radii"
[in print].

\item{7. }I. Ya. Bakelman, A. L. Verner, B. E. Kantor,
``Differential Geometry in the Large"[in Russian],
Nauka, Moskow, 1973.
       
\item{8. }W. Blaschke, ``Kreis und Kugel"
(Veit, Leipzig), 2nd Ed. De Gruyter, Berlin, 1956.

\item{9. } A. V. Pogorelov, ``Exterior Geometry of
Convex Surfaces" [in Russian], Nauka, Moscow, 1969.

\item{10. }W. Wiel, R. Schneider, ``Zonoids and
related Topics", in Convexity and its Applications,
Ed. P. Gruber and J. Wills, Birkhauser, Basel, 1983.

\hfill Institute of Mathematics

\hfill  Armenian Academy of Sciences

\hfill e.mail: rafik\@instmath.sci.am

\bye